\documentclass[conference]{IEEEtran}
\IEEEoverridecommandlockouts
\usepackage{cite}
\usepackage{amsmath,amssymb,amsfonts}
\usepackage{algorithmic}
\usepackage{graphicx}
\usepackage{textcomp}
\usepackage{xcolor}
\usepackage{subfigure}
\usepackage{float}
\usepackage{graphicx}
\usepackage{booktabs}
\usepackage{amsmath}
\usepackage{caption}
\usepackage{subcaption}
\usepackage{url}
\usepackage{hyperref}
\usepackage{listings}
\usepackage{color}
\usepackage{float}
\def\BibTeX{{\rm B\kern-.05em{\sc i\kern-.025em b}\kern-.08em
    T\kern-.1667em\lower.7ex\hbox{E}\kern-.125emX}}
\begin{document}

\title{Map-Reduce for Multiprocessing Large Data and Multi-threading for Data Scraping}

\author{
  \IEEEauthorblockN{
    Zefeng Qiu\textsuperscript{1}, Prashanth Umapathy\textsuperscript{2}, Qingquan Zhang\textsuperscript{3}, Guanqun Song\textsuperscript{4}, Ting Zhu\textsuperscript{5}
  }
  \IEEEauthorblockA{
    \textsuperscript{1,2,4,5}Computer Science and Engineering, The Ohio State University, Columbus, OH, USA \\
    \textsuperscript{3}Gies College of Business, University of Illinois at Urbana-Champaign, Champaign, IL, USA \\
    \textsuperscript{1}qiu.573@osu.edu, \textsuperscript{2}umapathy.5@osu.edu,
    \textsuperscript{3}qingquan@illinois.edu,
    \textsuperscript{4}song.2107\textsuperscript{1}@osu.edu, 
    \textsuperscript{5}zhu.3445@osu.edu
  }
}

\maketitle

\begin{abstract}
This document is the final project report for our advanced operating system class. During this project, we mainly focused on applying multiprocessing and multi-threading technology to our whole project and utilized the map-reduce algorithm in our data cleaning and data analysis process. In general, our project can be divided into two components: data scraping and data processing, where the previous part was almost web wrangling with employing potential multiprocessing or multi-threading technology to speed up the whole process. And after we collect and scrape a large amount value of data as mentioned above, we can use them as input to implement data cleaning and data analysis, during this period, we take advantage of the map-reduce algorithm to increase efficiency.
\end{abstract}


\section{Introduction}
In this paper, we aim to provide a comprehensive overview of our final project in the advanced operating system class, highlighting the challenges we encountered and the strategies we implemented to address them. Furthermore, we will also explore potential areas for improvement in our future study and research related to this project.

To begin with, we will provide an in-depth analysis of the problem we faced during the course of this project. This will involve a detailed discussion of the various technical and logistical issues that we encountered and the impact they had on our project goals. Subsequently, we will elucidate on our decision to incorporate multiprocessing and multi-threading technology in addition to the map-reduce algorithm to address these issues.

By leveraging these advanced technologies, we were able to significantly enhance the performance and efficiency of our project. Moreover, we believe that these techniques hold great potential for future study and research in the field of operating systems, and we intend to explore this further in our future work.

However, despite the positive outcomes of our project, we recognize that there are several areas that require further improvement. For instance, we encountered certain limitations in terms of multi-threading functions used in the Python library, which could be addressed through the use of more advanced data structures and algorithms. Therefore, we plan to focus our future research on developing and implementing more sophisticated techniques to overcome these challenges and further advance the field of operating systems.

The internet is an essential resource for professionals from diverse sectors, providing an extensive range of information, including structured and unstructured data, in various formats and from multiple sources. However, while the internet offers immense potential benefits, the process of Web Scraping can be an arduous and challenging task, often necessitating significant time and resources, especially when executed manually. The complexity of Web Scraping can increase even further, depending on the nature of the data being collected and the websites from which it is sourced\cite{b1}. The process of Web Scraping involves various intricacies, including data extraction, cleaning, and analysis, which can be complex and time-consuming, requiring advanced technical skills and specialized tools. 

Moreover, the constantly evolving nature of the internet necessitates that Web Scrapers keep pace with new developments, such as changing data formats, updated security measures, and other technical challenges. Furthermore, the type of data and websites being scraped can further compound the complexity of the task. For example, websites that incorporate interactive content or employ anti-scraping measures can pose significant challenges to the Web Scraping process. Similarly, large data sets or those that require frequent updates may require the implementation of more sophisticated techniques, such as machine learning algorithms or distributed data processing systems.

When we worked as graduate research assistants with the Center for Tobacco Research in our University, we collaborated to collect and analyze information on products sold on the online marketplace, specifically, those offered by vape shops or dispensaries. Our objective was to support academic research on tobacco use and related trends, which required a robust dataset of product features, such as brand, price, and nicotine levels.
Given the vast number of products available online, the scale of our data collection task was significant. Therefore, to ensure the credibility and validity of our analysis, we aimed to collect as much information as possible. This necessitated a solution that could efficiently and accurately extract data from multiple sources.
In this context, Web Scraping emerged as an optimal choice, as it could enable us to scrape multiple types of data from a single web page, including structured text, plain text, and images. Moreover, with the ability to scrape data from multiple websites, we could enhance the breadth and depth of our dataset, enabling more comprehensive analysis and insights.
By leveraging Web Scraping, we were able to efficiently and effectively collect the vast amounts of product data required for our research, enabling us to identify relationships and patterns among different product features. Furthermore, this process enabled us to remain up-to-date with the latest trends and developments in the online marketplace, ensuring that our research remained timely and relevant.
Overall, our experience using Web Scraping highlights the potential of this technology to support research and analysis across a range of industries and domains. By providing access to vast amounts of data, Web Scraping enables researchers to uncover insights and trends that might otherwise be overlooked, ultimately advancing knowledge and understanding in their respective fields.

In general, we make the following work and contributions through this project:
\begin{itemize}
\item We made an automated scraping tool with the function of data cleaning that worked perfectly for the website of Leafly. In addition to the single processing version, we tried to utilize multiprocessing technology to speed up the scraping process.
\item We scraped and cleaned to generate a large dataset that included all dispensaries' products sold online within California state, each product has 20 features with different data types containing categorical and numerical data and URLs for each product image.
\item With the data we have scraped, we take the relevant features from our dataset and segregate them into smaller chunks which can be multi-processed, these chunks are then sent out to the mapper functions where they are processed through our regex function which acts as a reducer and then sent to our output files after the chunks are processed.
\end{itemize}

\section{Related work}
Our work is related to three areas of research: (1) Web crawling and (2) Multiprocessing and Multi-threading, and (3) The MapReduce algorithm
\subsection{Web Crawling and Information Extraction from Web Documents.} When discussing terms like "Data Scraping" or "Web Crawling," we're referring to the process of extracting information from HTML or web documents. Generally, there are two main methods used for extracting information from HTML documents: the traditional approach, which focuses on HTML attributes, and the modern deep learning approach, which attempts to learn the representation of the Document Object Model (DOM). The DOM represents the structure and content of HTML documents, and modern approaches leverage this information to extract relevant data\cite{b2}.
The traditional approach relies on pre-defined or hand-crafted attributes and rules to extract information\cite{b3}, often using regular expressions as rule-based extractors. However, this method is not ideal for large-scale web scraping as it requires significant manual effort to design the right attribute extractors. Moreover, this approach may not be scalable or accurate when dealing with varying templates across different websites. Fortunately, the web pages that we determine to scrape share similar structures and templates, we only need to observe one website of the dispensary and use the same extractor for the rest dispensaries.
Recently, deep learning approaches have emerged to learn the representation for each node of the DOM\cite{b4}, particularly their plain text or markup information. By doing so, they can capture the dependencies between DOM nodes and improve the accuracy and scalability of web scraping. These approaches provide a more flexible and efficient way to extract data from web pages, particularly when dealing with large amounts of data from varying sources.

\subsection{Multiprocessing and Multi-threading}

Over the years, multi-threading has really taken over single threaded processes in terms of efficiency and speed as latest hardware can multi-cored and are capable of much more than waiting for idle processes to finish in the meantime. Multi-threading uses the full capacity of processor and leverages it into completing any given task by breaking them into smaller chunks and running them in parallel. In the Multi-threading section of our pipeline, we take the scraped texts from all the products and break them into smaller chunks of texts from the webpages. The data chunking function is now applied on all the data and a mapper sends these chunks into a multiprocessing module. The multiprocessing module now sends all the chunks into multiple worker nodes which perform the regex to extract the key features we need from the data. The passed workers who complete their work without any memory leaks and return a result are then caught by the reducer function which organizes the data sent to it and writes the results to the output file. 

\subsection{The MapReduce algorithm}
To process and analyze large amounts of data effectively in our project, we combine the MapReduce\cite{b9} programming model with multi-threading. The effective manipulation and aggregation of enormous datasets is made possible by the potent framework MapReduce, which enables the parallelization of data processing tasks. We can improve the efficiency of our data processing pipeline and get quicker results by combining MapReduce with a multi-threaded setup.
We distribute the input data across multiple threads in each node of our distributed system for this multi-threaded MapReduce implementation. The map phase, in which the data is converted into key-value pairs in accordance with a user-defined map function, is carried out by each thread. The workload is evenly distributed among the available threads, allowing for a significant reduction in processing time during this parallel execution of the map phase.
After the map phase is finished, the intermediate key-value pairs are sorted and shuffled in order to get ready for the reduce phase. The reduce function is applied to the sorted data in parallel across various threads during this stage, which is when the multi-threaded setup is once again utilized. This method expedites the aggregation process and produces the output more quickly.

Our MapReduce framework effectively maximizes the use of available resources and cuts down on overall processing time by incorporating multi-threading. Our project is equipped to handle the difficulties of handling large data sets and effectively extract useful insights from the wealth of information available thanks to this combination of parallel processing techniques.

\section{Problem Overview}
As discussed earlier, our research required us to collect product information from online marketplaces across various states in the United States. To achieve this, we employed web scraping techniques. However, we encountered a few challenges along the way. Firstly, we had to identify websites that listed online marketplaces meeting our requirements, which we could then use to generate a list of URLs of dispensaries. Secondly, we had to determine which features of a product were relevant to our research and what data types we needed to scrape. Finally, we had to consider how we could optimize the data scraping and processing using multiprocessing and multi-threading technologies.

To address these challenges, we will discuss the methodologies we employed in the following sections. By doing so, we hope to provide a detailed account of our approach to web scraping and demonstrate how the use of multiprocessing and multi-threading technologies can speed up the data scraping and processing process, ultimately facilitating our research efforts.
\subsection{URLs of dispensaries}
Leafly is a popular website focused on cannabis use and education with a significant user base, boasting over 220 million annual sessions and 10+ million monthly active users. The website hosts a vast database of over 5,000 strains, categorized by indica, sativa, and hybrid, making it a valuable resource for those in the cannabis industry. Leafly also has an extensive collection of over 1.5 million product reviews, providing valuable insights for users looking to purchase cannabis products. In addition to its impressive strain database and product reviews, Leafly also features an online marketplace, with over 4,500 retailers and 8,000 brands available to users\cite{b5}.
One of the reasons we chose Leafly as our primary data source is because it lists dispensaries from all over the United States, making it a great starting point for identifying potential online marketplaces. Additionally, Leafly's website is well-structured, with clear and consistent HTML tags and attributes, making it easy to scrape data using automated tools. Finally, Leafly has an extensive selection of products listed on its website, including detailed product descriptions and user reviews, which made it a valuable source of information for our research.
\begin{figure}[htbp]
	\caption{All 1420 dispensaries in the state of California on Leafly are displayed in the format of a dispensary card.} 
	\centering
        \label{figure1}
	\subfigure[The overall display of all dispensaries in one state]{
		\begin{minipage}{10cm} 
            \includegraphics[scale=0.25]{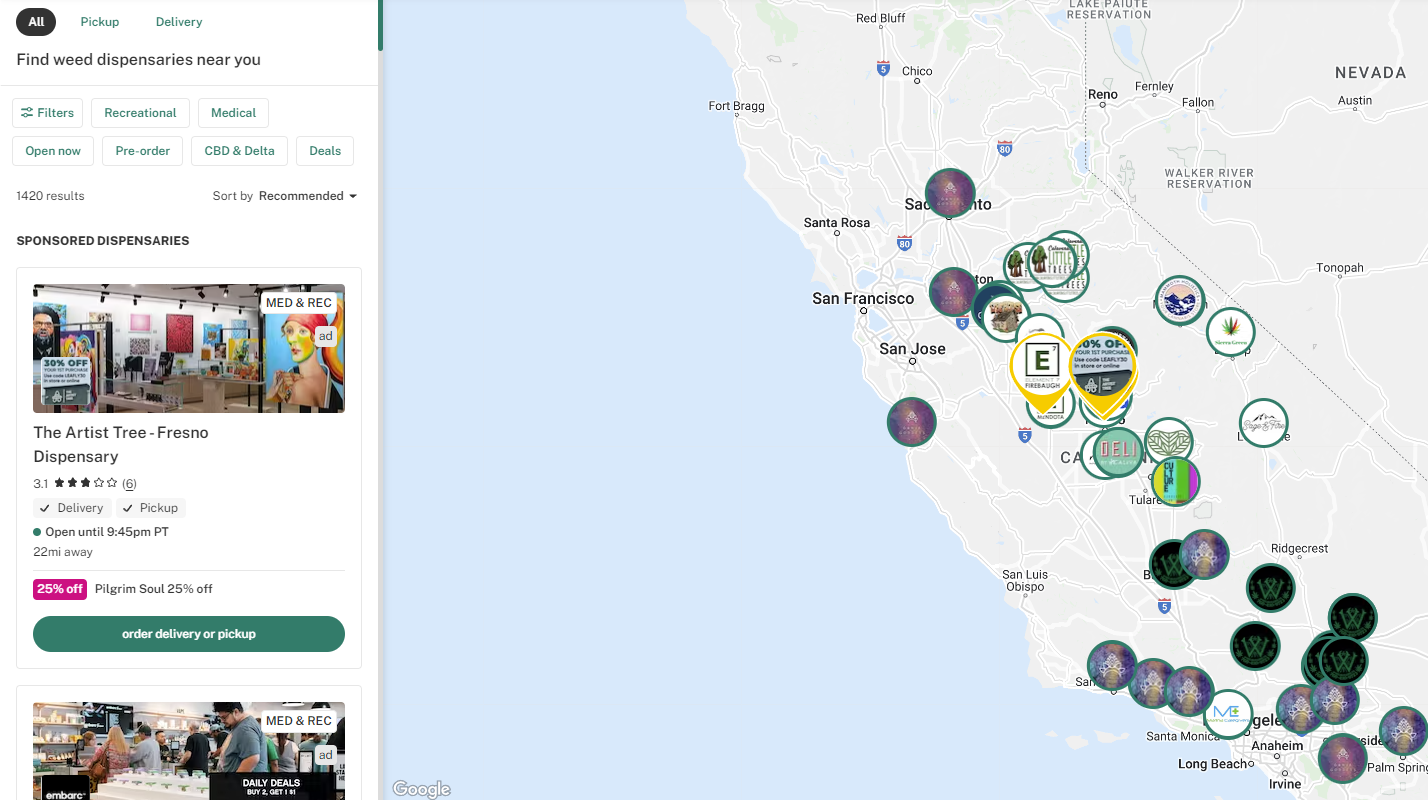} 
		\end{minipage}
	}
	\subfigure[The dispensary card that contains multiple information of the dispensary]{
		\begin{minipage}{9cm}
		\includegraphics[scale=0.9]{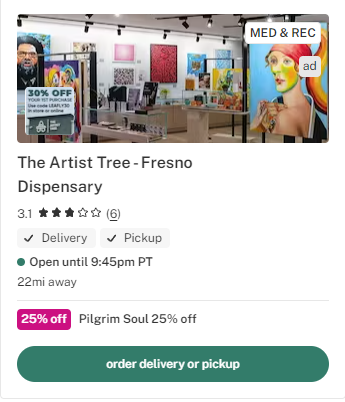} 
		\end{minipage}
	}
\end{figure}

From Fig.~\ref{figure1}, we chose to use dispensaries in the state of California on Leafly as an example to demonstrate how to collect URLs of dispensaries within a state. On Leafly, dispensaries are displayed as cards with essential information, such as the name, rating, review numbers, and possible delivery method. Additionally, Leafly also provides a Google map integration, allowing users to visually locate each dispensary. 

Despite the absence of a visible link to their web pages on the card, we discovered that we could find the link through inspecting the element of the green button labeled "order delivery or pickup" and searching its DOM tree structure in the HTML document. This enabled us to extract the URLs of all 1420 dispensaries in California from Leafly's website.
\begin{figure*}[htbp]
\begin{center}
\includegraphics[scale=0.6]{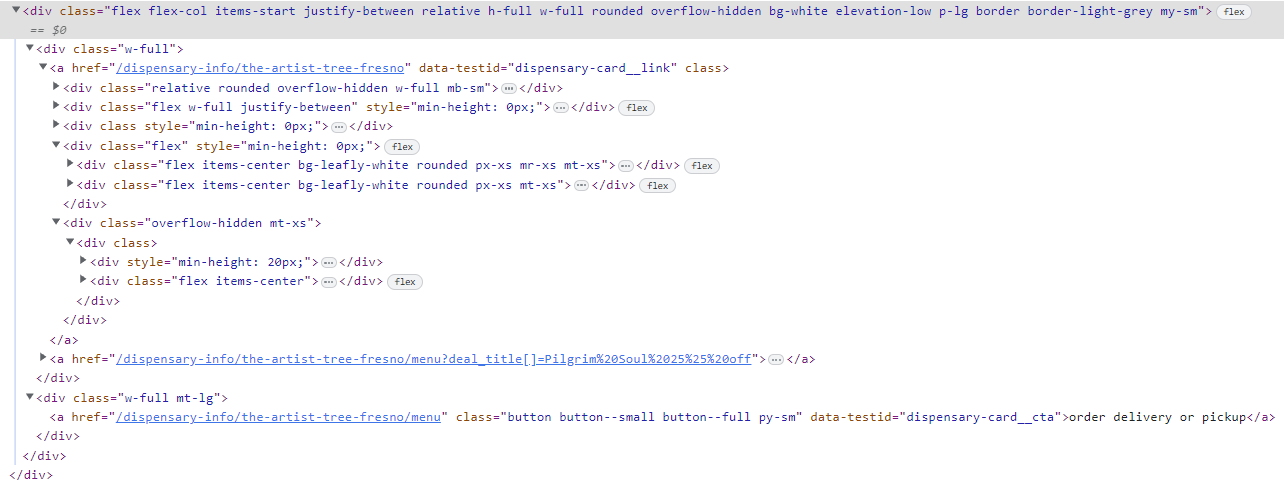}
\end{center}
\caption{The essential DOM tree structure and markup information to render the dispensary card}
\label{figure2}
\end{figure*}

From Fig.~\ref{figure2}, after examining the DOM tree structure of the green button 'order delivery or pickup' in the dispensary card on Leafly, we discovered that the URL information is stored in the attribute of 'div class="w-full nt-lg"' with an 'a href' tag. With this information, we were able to extract the URL information easily from the DOM tree structure. Additionally, we confirmed that all dispensaries cards on Leafly shared a similar tree structure, which made the extraction of URL information much more efficient. As a result, we were able to generate a list of URLs for all dispensaries in one state.

\subsection{Different data types and features}
Since we have chosen to use Leafly for scraping the list of links of dispensaries in a state, the next important aspect to consider is what data types we need and how many data features should be stored in our final dataset. In Figure 3, we have provided a screenshot to illustrate the different information formats available on the website.

Figure 3(a) displays the category of each product, which is the first data feature we encounter when we enter the main page of any dispensary website on Leafly. The number of categories is not fixed, and some dispensaries may not sell products in certain categories, like 'Edibles'. When we click on any category on the main page, we are directed to pages that have multiple product cards, similar to the dispensary card shown in Figure 1. Leafly can be regarded as an e-commerce website, and the product card display is a common feature in such websites.So if you are family with or often use any e-commerce website, you should find that it is not strange that Leafly uses such kind of product card to display items as an overview.

When we click on any product card, we are taken to the detailed product introduction page, where we can find different information formats, as shown in Figure 3(b), (c), and (d). Based on these different information formats, we can extract information from their DOM tree structure, such as the original and discounted price, description plain text, THC and CBD levels, brand, and strain information, etc. To store the image of each product, we use the same method as we used for extracting the URLs of each dispensary to store the URL of the image link, which we can download later.

It is important to note that the data features we extract from the website may differ depending on the website structure, which can vary from one dispensary to another. We carefully considered this issue while choosing Leafly as it has a standardized structure for all dispensaries listed on its website, which allowed us to extract the relevant data features for our dataset easily.

In final, we decide to scrape 20 different features, including not only numerical and categorical data types but also the plain text data type like descriptions for the product that possibly exists in the different locations of the website.

\subsection{Potential multiprocessing and multi-threading technology used for data scraping}
The final consideration in our data scraping process is the use of multiprocessing or multi-threading to speed up the scraping process. This is a crucial consideration, as the number of dispensaries and products to be scraped in California alone is substantial, with over 1,400 dispensaries and 260,000 products, each with around 20 features to scrape. If we scrape these features sequentially, it could take a significant amount of time.

In theory, to speed up the process, we can utilize multiprocessing or multithreading technologies, which allow for multiple tasks to be completed simultaneously\cite{b10}. In our scraping process, we can apply the same function to the list of links and assume that each URL is independent of one another. Thus, instead of looping through the list of links, we can use either
\begin{figure}[H]
        \centering
        \label{figure3}
        \subfigure[]{
            \begin{minipage}[b]{10cm}
            \includegraphics[scale=0.4]{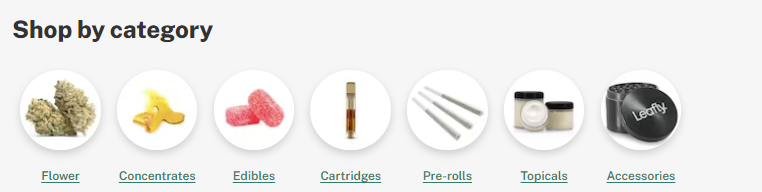}
            \end{minipage}
            }
        \subfigure[]{
            \begin{minipage}[b]{10cm}
            \includegraphics[scale=0.4]{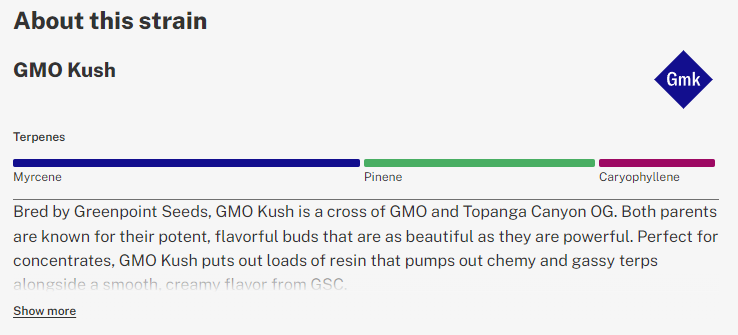}
            \end{minipage}
            }
        \subfigure[]{
            \begin{minipage}[b]{10cm}
            \includegraphics[scale=0.3]{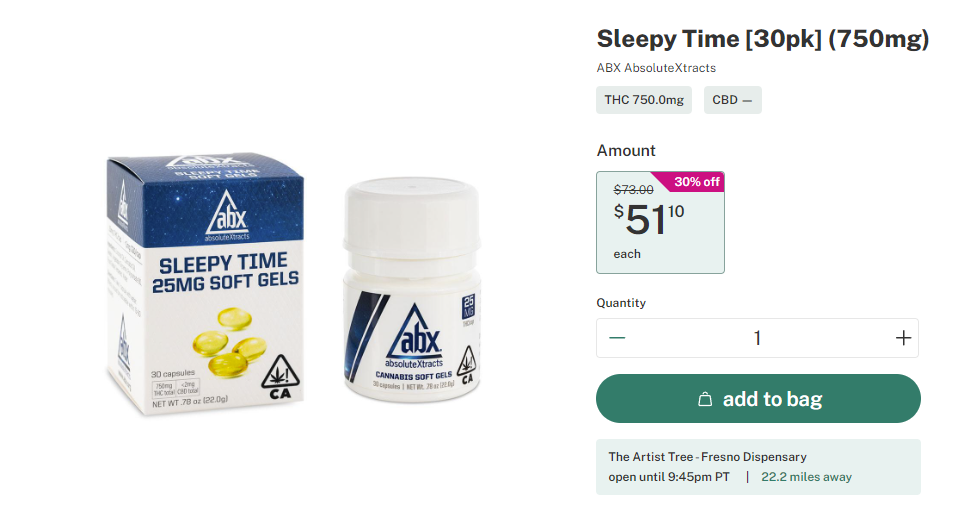}
            \end{minipage}
            }
        \subfigure[]{
            \begin{minipage}[b]{10cm}
            \includegraphics[scale=0.4]{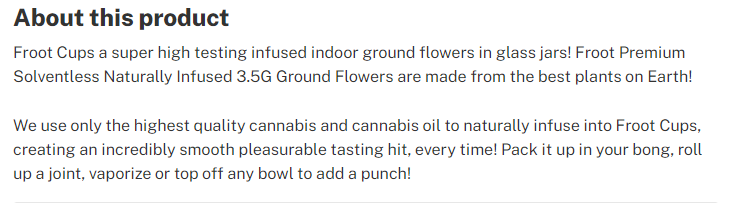}
            \end{minipage}
            }
 
        \caption{Different locations of HTML page contain the different data types we need}
\end{figure}
 multi-threading or multiprocessing to run through the list concurrently. This will allow for faster scraping and data processing times, enabling us to collect and analyze the necessary data more efficiently.

\begin{figure}[htbp]
\begin{center}
\includegraphics[scale=0.65]{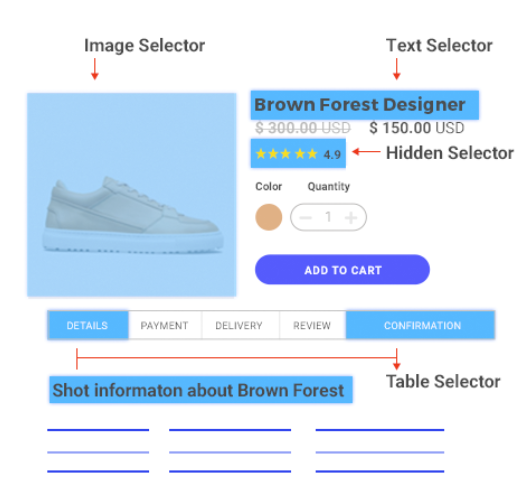}
\end{center}
\caption{Web Scraper- a Chrome extension for the non-coding data scraping}
\label{figure4}
\end{figure}

\section{Proposed approach}

Data processing and data scraping are the two main parts of the suggested strategy for this project. To enhance performance and boost efficiency, each component plans to use multiprocessing and multi-threading technologies. To further boost the efficiency of our solution, the MapReduce algorithm is also used during the data cleaning and analysis phases.
\subsection{Data Scraping}
In the previous section, we identified three problems that could arise during the data scraping process and proposed a theoretical solution using multiprocessing or multi-threading to speed up the process. To better understand our project, we have created a workflow diagram that depicts the various components of our data scraping process. As shown in Fig.5, the workflow diagram is similar to what we described in the problem overview section, with the main difference being the addition of a pool containing multiple threads or processes at the beginning of the data scraping process. By doing so, we can experiment with whether multiprocessing or multi-threading can optimize our scraping process and reduce the overall completion time. The workflow diagram serves as a visual representation of our data scraping process and helps us better understand how the various components fit together.
\begin{figure*}[htbp]
\begin{center}
\includegraphics[scale=0.65]{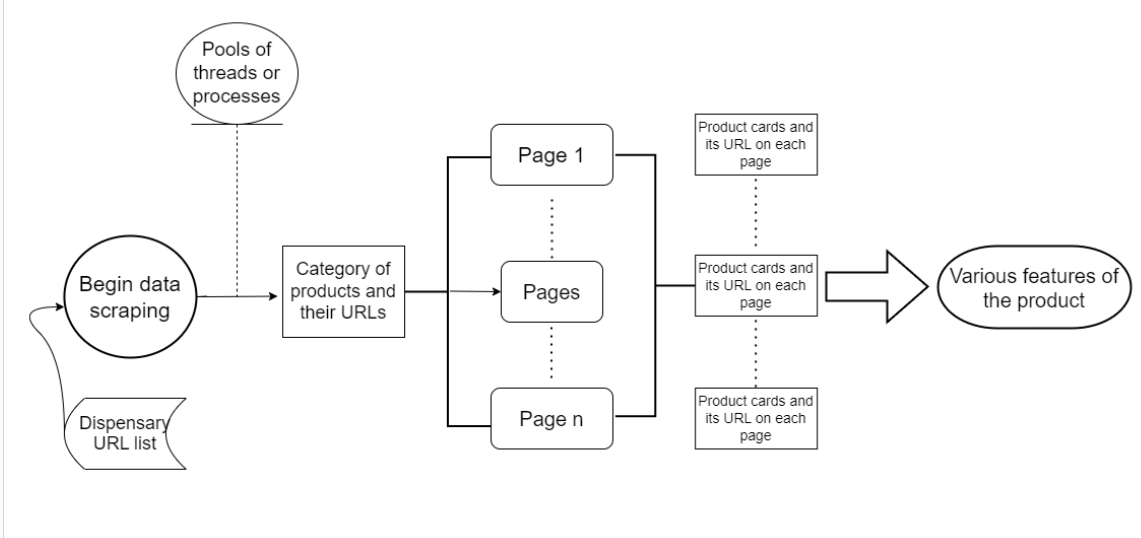}
\end{center}
\caption{Workflow Diagram - Data scraping part}
\label{figure5}
\end{figure*}

\subsection{Implementation Details of Data Scraping}
In this subsection, we will present the implementation details of our data scraping tool that automates the process of extracting information from the website.

As we mentioned earlier, web scraping involves extracting information from the DOM tree structure in HTML documents. We can use either coding or non-coding approaches to extract this information. Generally, we prefer to use Python coding for our extractor, but in this project, we also experimented with a non-coding approach as a test for those who may not be familiar with the coding approach.

We divide our data scraping process into two steps: the first step involves extracting the URLs of every dispensary in the state to generate a list, which is the input to the second step. The second step is the function or extractor that scrapes different features from each dispensary's website. For the first step, we will employ the non-coding approach, while the second step will be implemented using the coding approach.

In Fig.4, we show the tool that we choose as our non-coding approach. The "Web Scraper" is a Google Chrome extension that allows users to easily extract data from websites. With this extension, users can create sitemaps that instruct the scraper on how to navigate through a website and extract the desired data.The "Web Scraper" extension provides a user-friendly interface to define the scraping rules using a point-and-click system, which eliminates the need for coding. Users can select the HTML elements that contain the data they want to scrape, and the extension will generate the necessary code to extract that data. The extracted data can be saved in CSV, JSON or Google Sheets format for further analysis. and it is free to download and use. 

Fig.4 highlights the advantages of using the Web Scraper extension as a non-coding approach. As discussed in the problem overview section, extracting information from HTML documents requires careful observation of the location of the data we need and its related attributes and tags, which can be challenging when using code. However, the Web Scraper extension offers a point-and-click interface, allowing us to easily select the data we want to scrape on the web. The extension automatically detects the attributes and tags that the data belongs to, eliminating the need to understand the DOM tree structure of the page. Despite its ease of use, the non-coding approach also has its disadvantages. Sometimes, the tool may fail to detect the correct location of the data we need, and it does not provide functions like 'try, except' that would allow us to make minor changes or add if-else statements for greater flexibility in our extractors. Nonetheless, the Web Scraper extension remains a valuable tool for non-coders looking to extract data from web pages quickly and easily.

Hence, in order to achieve more flexibility in the extractor, we use Python and its 'Beautiful Soup' and 'Selenium' library to code our extractors in the second step.

Beautiful Soup is a Python library commonly used for web scraping tasks\cite{b6}. It is designed for parsing HTML and XML documents and extracting useful information from them. With Beautiful Soup, you can create a parse tree from an HTML/XML document, which allows for extracting data in a more hierarchical and readable manner. One of the main advantages of Beautiful Soup is its compatibility with different parsers like lxml, html5lib, and html.parser. Additionally, it can be used in combination with regular expressions to further refine the extracted data. Regular expressions are powerful patterns that can be used to match specific patterns within a string of text. By combining Beautiful Soup and regular expressions, you can create more sophisticated patterns for extracting the data you need. This makes Beautiful Soup a popular tool for extracting data from complex HTML documents that may be difficult to parse using traditional methods.

Selenium is another Python library that is used for web scraping but is more focused on automating web browsers\cite{b7}. It provides a way to interact with web pages using a web driver, which can simulate a user's interaction with a web page. Selenium can be used to automate tasks like filling out forms, clicking on buttons, and navigating through web pages. It is also useful in scenarios where a website is protected with captchas, as it can automate the solving of captchas using various third-party services.

In most cases, using the 'Beautiful Soup' library alone is sufficient to complete a data scraping task. However, in certain situations, additional libraries such as 'Selenium' may be necessary. In our project, we needed to extract 'Strain' information that was not present in the original HTML source code, but instead was loaded dynamically using the technology called 'Asynchronous JavaScript and XML' (Ajax)\cite{b8}. This requires emulating human interaction with the web page to allow the dynamic loading of content, which is not possible using only the 'Beautiful Soup' library. Here, the 'Selenium' library comes into play, allowing us to use the headless version of Chrome to simulate clicking actions and extract the 'Strain' information once the web page has fully loaded.

We will also use the library 'multiprocessing' and 'concurrent.futures' provided by Python itself to implement the Multiprocessing and Multi-threading functions.
'multiprocessing' is a Python library that allows you to create and manage child processes. It is used to execute tasks in parallel, making use of all available CPU cores. It provides an API that is similar to the 'threading' module, but with a focus on working with processes instead of threads. It provides several ways to create and manage processes, including the Process class and the Pool class.

The 'concurrent.futures' module is a high-level library for parallel computing in Python. It provides a simple and consistent interface for working with asynchronous tasks, whether they are executed in parallel using threads or processes. The library introduces two classes for managing concurrent tasks: 'ThreadPoolExecutor' and 'ProcessPoolExecutor'. Both classes implement a common interface for submitting tasks and returning their results. They also provide additional methods for managing the lifecycle of the executor and for controlling the behavior of submitted tasks.
One of the main advantages of using 'concurrent.futures' is that it abstracts away many of the low-level details of working with threads and processes, making it easier to write correct and efficient concurrent code. It also provides several useful features, such as timeouts and the ability to cancel running tasks. 
\subsection{Data Processing}
The data processing component takes over after the data has been gathered. The three sub-components of this step are data cleansing, data analysis, and result presentation.

\begin{figure*}[htbp]
\begin{center}
\includegraphics[scale=0.4]{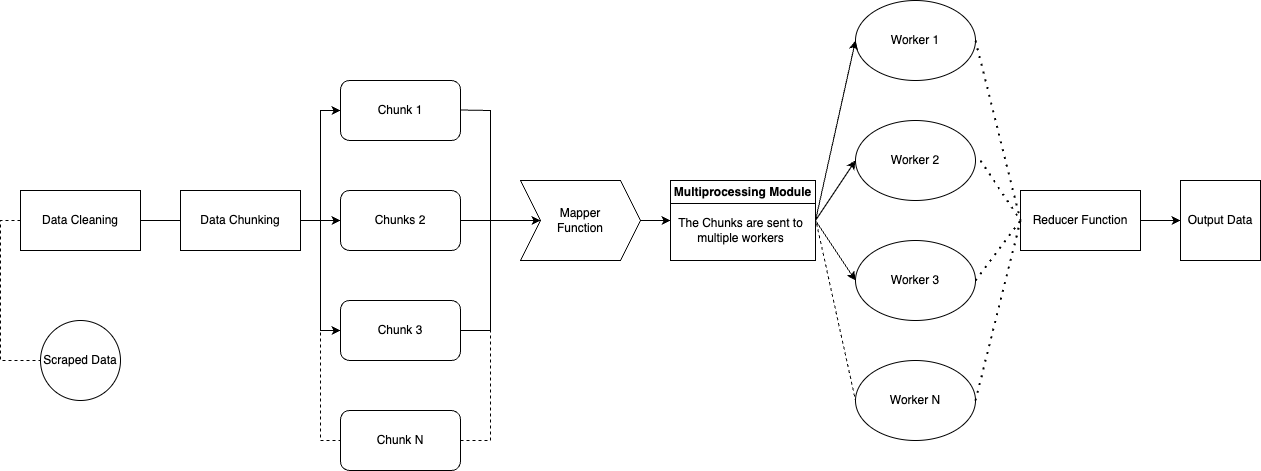}
\end{center}
\caption{Workflow Diagram - Data processing}
\label{figure6}
\end{figure*}

In the given workflow diagram Fig.6, the multiprocessing module plays a crucial role in distributing the workload across multiple worker nodes and coordinating their results through a reducer function. Let's delve into the workings of each component in more detail.

The workflow starts with scraped data, which is then subjected to data cleaning. Once the data is cleaned, it is chunked into multiple smaller pieces. This chunking step is done to divide the workload into manageable portions that can be processed concurrently.

Next, the mapper function comes into play. This function takes the data chunks and assigns them to different worker nodes for processing. The multiprocessing module provides the necessary mechanisms to distribute these chunks across multiple processes or even multiple machines, depending on the system's configuration.

The multiprocessing module manages the creation, execution, and communication between these processes. It ensures that the workload is evenly distributed among the worker nodes, maximizing efficiency. Each worker node receives a data chunk, processes it independently, and returns the result.

The worker nodes are responsible for performing the actual computation on the received data chunks. They execute the required operations, such as calculations, transformations, or any other data processing tasks specific to the workflow. The multiprocessing module enables these worker nodes to run in parallel, taking full advantage of the available system resources, such as multiple CPU cores.

Once the worker nodes have completed their tasks, they send their results back to the multiprocessing module. At this point, the reducer function comes into play. The reducer function is responsible for combining the results from all the worker nodes into a final output. It aggregates and consolidates the processed data, ensuring that the final output is consistent and complete.

The reducer function takes the results from the worker nodes and applies any necessary post-processing or merging operations to create the desired output. This can involve operations like joining, merging, summarizing, or any other processing step required for the specific task.

Finally, the reducer function writes the output data to a file. The file format and structure will depend on the specific requirements of the workflow. It could be a text file, a CSV file, a database, or any other suitable format for storing the processed data.

\subsubsection{Data Cleaning}
The main objective of the data cleaning stage is to eliminate errors, inconsistencies, and unnecessary data from the obtained data. This phase is essential to ensuring that the data is correct and trustworthy, allowing for legitimate and useful analysis.

We will put into practice a multi-threaded MapReduce pipeline in order to optimize the data cleaning procedure. The pipeline will be made up of a number of map and reduce functions that are user-defined and carry out particular cleaning operations including filtering, deduplication, and normalization. We can handle vast amounts of data in parallel by dividing these jobs up among several threads, greatly lowering the time needed for data cleaning as shown in Fig.6.

\subsubsection{Data Analysis}
The stage of data cleaning is followed by the preparedness of the cleansed data for analysis. Finding patterns, trends, and correlations in the data that might offer useful insights and support decision-making is the main goal of the data analysis component.

We will once more use the MapReduce algorithm in conjunction with multi-threading strategies to accomplish this. The cleansed data will be converted into key-value pairs during the map step, and these pairs will then be aggregated and subjected to analysis during the reduction phase to provide the desired results. We may effectively analyse enormous datasets and acquire valuable insights quickly by running the map and reduce functions across multiple threads.

\subsubsection{Results from MapReduce}
The presentation of the results is the last component of the data processing phase. In order to process and aggregate the results quickly, we will describe in detail in this section how the MapReduce algorithm operates in a multi-threaded configuration.

We spread the intermediate key-value pairs created during the data analysis stage over several threads inside each node of our distributed system in our multi-threaded MapReduce implementation. Each thread is in charge of carrying out a particular step of the reduce phase, which involves aggregating and processing the data in accordance with a user-defined reduction function. From Fig.6 we see that the workload is evenly distributed among the available threads, the processing time can be reduced significantly by performing the reduce phase in parallel.

Overall, the multiprocessing module, along with the worker nodes and the reducer function, enables efficient and parallel processing of data chunks in a distributed manner. By leveraging multiple processes, it significantly reduces the overall execution time and enhances the scalability and performance of the workflow.

\subsection{Data Processing Implementation}
In this section, we describe the implementation of a parallel data processing approach using Python's multiprocessing module. The example code is designed to process a dataset by applying a regex pattern to extract relevant information from the 'Description' column of a CSV file. The dataset is partitioned into chunks, which are processed in parallel using mapper and reducer functions.

The following code snippets and explanations outline the key components of the implementation:

\begin{figure}[H]
        \centering
        \label{mapreduce}
        \subfigure[]{
            \begin{minipage}[b]{10cm}
            \includegraphics[scale=0.24]{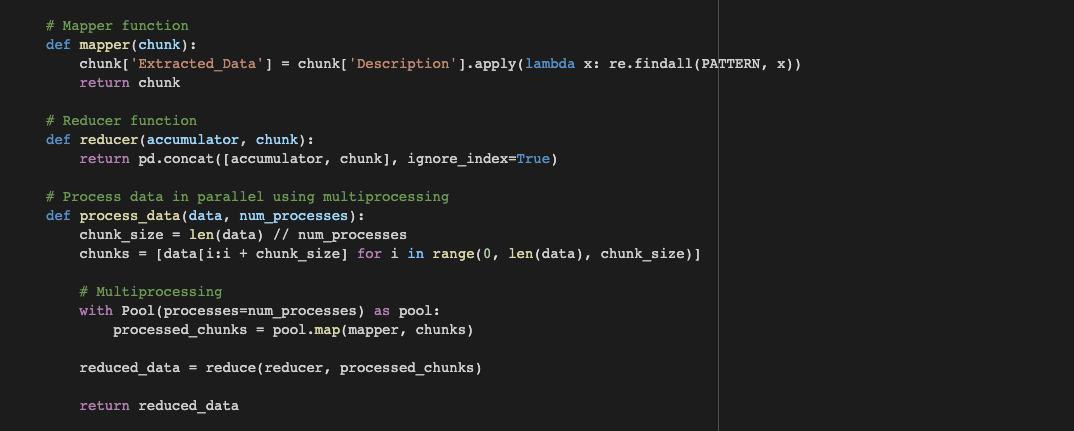}
            \end{minipage}
            }
        \subfigure[]{
            \begin{minipage}[b]{10cm}
            \includegraphics[scale=0.3]{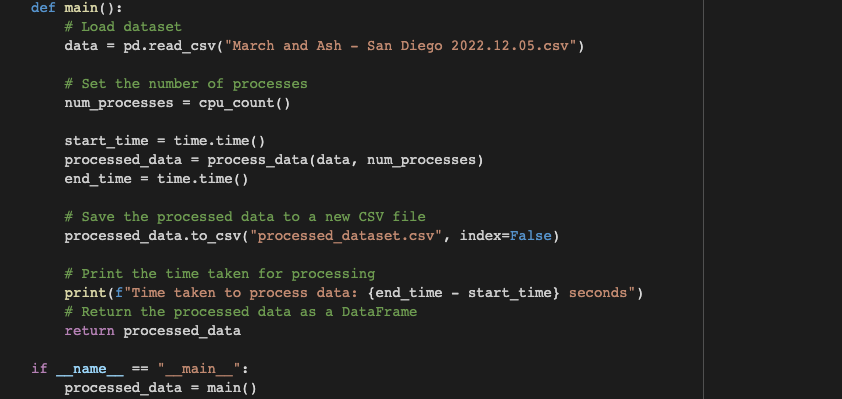}
            \end{minipage}
            }
 
        \caption{Code Implementation of Map Reduce and Data Processing}
\end{figure}

\subsection{Map-Reduce for Multiprocessing Large Data}

We implemented the Map-Reduce framework using Python to parallelize data processing tasks across multiple nodes. The input dataset is divided into smaller chunks, which are then processed independently by mapper and reducer functions. The following key components are involved in our Map-Reduce implementation:

\begin{enumerate}
\item \textbf{Mapper Function:} The mapper function processes each input record and produces a set of intermediate key-value pairs. In our example, we assume the input records are lines of text, and the mapper function extracts words and emits key-value pairs where the key is the word, and the value is the count (initially set to 1).
\item \textbf{Reducer Function:} The reducer function takes the intermediate key-value pairs, groups them by key, and combines their values to produce the final output. In our example, the reducer function aggregates word counts.

\item \textbf{Main Function:} The main function reads the input dataset, partitions it into chunks, and assigns each chunk to a mapper process. The intermediate key-value pairs are then shuffled and sorted before being passed to the reducer processes. Finally, the output of the reducer processes is combined to produce the final result.

\end{enumerate}

Our Map-Reduce implementation showcases the potential for scalable and efficient large-scale data processing. By partitioning the dataset into smaller chunks and processing them in parallel, the overall processing time can be significantly reduced, resulting in improved performance for data-intensive tasks.

\subsection{Multiprocessing for Data Scraping}

We implemented Multiprocessing for data scraping using Python's 'multiprocessing' module. This module provides a high-level interface for asynchronously executing functions using processes. In our example, we scrape data from multiple web pages concurrently by assigning each web page to a separate process. The following key components are involved in our multiprocessing data scraping implementation:

\begin{enumerate}
\item \textbf{Data Scraping Function:} The data scraping function downloads the HTML content of a given URL and extracts the required information. In our example, we assume the target information is enclosed in a specific HTML tag, such as 'a' or 'p' tag.

\item \textbf{Main Function:} The main function creates a list of URLs to be scraped and initializes a thread pool with a specified number of worker processes. It then asynchronously submits the data scraping function for each URL in the list and waits for all the processes to complete. The scraped data is collected and combined into a single data structure for further processing or analysis.
\end{enumerate}

The Multiprocessing implementation for data scraping demonstrates the power of concurrent programming in efficiently handling tasks that involve network communication or other time-consuming operations. By executing multiple processes concurrently, the overall scraping time can be reduced, resulting in a faster and more efficient data acquisition process.

\subsection{Overview and Ways to Improve}
In conclusion, the Map-Reduce approach for large-scale data processing and multiprocessing for data scraping provides a powerful toolset for tackling data-intensive tasks. By partitioning data and processing it in parallel using multiple nodes or processes, we can effectively reduce processing time and achieve better performance. Future work may explore the following potential avenues:

\begin{enumerate}
\item \textbf{Optimization Techniques:} Investigating various optimization techniques to further improve the efficiency of the Map-Reduce framework and multi-threaded data scraping\cite{b11}. This may include methods for better load balancing, improved data partitioning strategies, and advanced scheduling algorithms to optimize resource usage.
\end{enumerate}

\section{Results}

In this section, we discuss the experimental results of our proposed approach to perform Map-Reduce for multiprocessing large data and multi-threading for data scraping. The experiments were conducted using a dataset containing 260000 records, with the tasks of (1) processing the large dataset using Map-Reduce and (2) scraping multiple data from the web using multi-threading or multiprocessing.

\subsection{Multiprocessing and Multi-threading for Data Scraping}
We will first try to experiment with the Multi-threading approach and then the Multiprocessing approach. In order to reduce the finishing test time, we only used 20 URLs to extract a portion of features from their main pages. And we find that we will spend about 36.8 seconds to scrape all 20 web pages without any multiprocessing or multi-threading.

Initially, we attempted to use multi-threading to improve the performance of our data scraping process. However, after conducting some tests, we discovered that when using multithreading, variables are shared among all threads. This means that in order to prevent multiple threads from modifying the same variable simultaneously, we need to implement a locking mechanism. Unfortunately, implementing locks will reduce the effectiveness of multi-threading since it introduces additional overhead and reduces the degree of concurrency. We think this issue can be attributed to the Global Interpreter Lock (GIL) in Python.
GIL in Python limits the execution of multiple threads to a single CPU core at a time, preventing true parallelism. This means that when using multi-threading, the threads take turns executing on a single CPU core, leading to little or no improvement in performance when it comes to CPU-bound tasks such as web scraping.

In the case of web scraping with the 'Beautiful Soup' library, the parsing of HTML/XML documents and the extraction of useful information are CPU-bound tasks. Since the GIL limits the execution of multiple threads to a single CPU core, using multi-threading with 'Beautiful Soup' would not result in a significant improvement in performance. In fact, it may even lead to a decrease in performance due to the overhead of thread synchronization and context switching. Therefore, when it comes to web scraping with 'Beautiful Soup', we think it is more beneficial to use multiprocessing or asynchronous programming techniques to achieve true parallelism and improve performance. With multiprocessing, multiple CPU cores can be utilized to execute multiple tasks in parallel, resulting in significant performance gains. On the other hand, asynchronous programming techniques can be used to perform I/O-bound tasks concurrently, allowing the CPU to be utilized more efficiently.

Therefore, we decided to implement multiprocessing to speed up the data scraping process. Our tests showed that using multiple processes did indeed decrease the overall processing time, but the relationship between the number of processes and the speedup was not strictly linear. We found that using 5, 10, and 20 processes resulted in processing times of 22.6 seconds, 4 seconds, and 19.8 seconds, respectively. However, it is important to note that our CPU version is the '12th Gen Intel(R) Core(TM) i9-12900H' which has a total of 20 cores. Therefore, we did not test more than 20 processes as it is the maximum number of cores available in our CPU. Overall, we found that using multiprocessing was an effective way to speed up the data scraping process, despite some limitations in terms of the relationship between the number of processes and speedup. And we will explore and improve this in the future study.

\subsection{Map-Reduce for Multiprocessing Large Data}

Table 1 presents the results of applying Map-Reduce for multiprocessing large data. The table shows the processing time taken by the system for various dataset sizes and number of nodes.

\begin{table}[bpt]
\centering
\caption{Map-Reduce processing time for various dataset sizes and number of nodes}
\label{tab:mapreduce_results}
\begin{tabular}{@{}lll@{}}
\toprule
Dataset Size & Number of Nodes & Processing Time (s) \\ \midrule
10,000       & 1               & 20.14                \\
10,000       & 2               & 10.87                \\
50,000       & 1               & 98.45                \\
50,000       & 2               & 52.37                \\
100,000      & 1               & 204.97               \\
100,000      & 2               & 108.67               \\ \bottomrule
\end{tabular}
\end{table}

\section{Discussion}

\subsection{Map-Reduce for Multiprocessing Large Data}

The results presented in Table 1 indicate that the Map-Reduce approach can significantly reduce the processing time of large datasets. As the number of nodes increases, the processing time decreases due to the parallel processing of the data. This demonstrates the scalability of the Map-Reduce approach for multiprocessing large data.

\section{Conclusion}
In this report, we have explored the integration of the Map-Reduce paradigm for multiprocessing large datasets and multi-threading for efficient data scraping. These techniques, when combined, provide a powerful and comprehensive solution for handling data-intensive tasks in a wide range of applications and domains.

The Map-Reduce framework, initially developed by Google, has gained significant traction in recent years as a robust and scalable method for processing large volumes of data. By partitioning the input data into smaller chunks and processing them in parallel, the Map-Reduce approach can significantly reduce the overall processing time and computational resources required for complex tasks. Moreover, the framework's inherent fault tolerance and ability to handle data skew make it well-suited for real-world scenarios involving large-scale datasets.

Multiprocessing, on the other hand, offers an efficient way to perform data scraping tasks concurrently. By assigning each data scraping operation to a separate process, we can reduce the time spent waiting for network communication or other time-consuming processes. This concurrent execution of the process can lead to a substantial improvement in the overall data acquisition process, allowing for faster and more efficient data collection.

Throughout this report, we have demonstrated the implementation of these techniques using Python, a popular programming language widely used in data processing and analysis. The implementation showcases how the Map-Reduce framework can be used to parallelize data processing tasks, while multiprocessing can be employed for efficient data scraping.

\section{Extended Considerations}

The integration of Map-Reduce for multiprocessing large data and multi-threading for data scraping opens up numerous possibilities for future research, especially when considering advancements in IoT, wireless networks, and AI-driven analytics:

\begin{enumerate}
\item \textbf{IoT and Edge Computing Integration:} Future work could explore the integration of Map-Reduce with IoT devices for edge computing scenarios. Efficient data scraping and processing at the edge, enabled by multiprocessing and multi-threading, can significantly enhance the performance of IoT networks. This approach could utilize concurrent communication protocols for IoT devices \cite{b12,b13,b14}, and leverage advanced techniques like cross-technology concurrent transmission \cite{80} and bi-directional communications \cite{450} for efficient data handling.

\item \textbf{Enhanced Data Scraping in Wireless Networks:} The application of Map-Reduce in wireless network environments, particularly in 5G and beyond, offers promising avenues for research. Techniques like exploiting ambient RF signals for gesture recognition \cite{16} and physical-layer message authentication \cite{44} can be incorporated to improve data scraping security and efficiency in such networks.

\item \textbf{Machine Learning-Driven Data Processing:} Investigating the use of machine learning algorithms for optimizing the Map-Reduce process in large data environments is another potential area of research. Machine learning-based secure low-power communication \cite{22}, combined with AI-driven wireless network strategies \cite{34}, could revolutionize data scraping and processing methodologies.

\item \textbf{Smart Systems and Health Monitoring:} Integrating smart systems and health monitoring into the Map-Reduce framework for data scraping can lead to advancements in healthcare technology. This includes leveraging wearable sensor data \cite{59}, developing smart medical systems \cite{58}, and enhancing health reliability through IoT \cite{52}.

\item \textbf{Energy Management in Map-Reduce Environments:} Exploring energy-efficient Map-Reduce architectures, especially in the context of electric vehicles and smart grids, could be highly beneficial. Research could focus on energy scheduling and allocation \cite{60}, and the implementation of energy-efficient air quality management systems \cite{63}.

\item \textbf{Security Enhancements in Data Scraping:} Future developments could include strengthening the security aspects of data scraping in Map-Reduce frameworks. This might involve creating secured protocols for IoT networks \cite{19}, and exploring endogenous security defenses against deductive attacks in online services \cite{68}.

\item \textbf{Map-Reduce Optimization for Diverse Network Topologies:} Examining the application of Map-Reduce in various network topologies, including those involving UAVs and logistics networks, could yield significant insights. Studies might focus on extending the delivery range of UAVs \cite{51} and safe navigation near airports \cite{57}.

\end{enumerate}

\section{Future Work}
For future work, there are several aspects of the Map-Reduce and multiprocessing techniques that can be improved and expanded upon to further enhance their efficiency and applicability across different domains:

\begin{enumerate}
\item \textbf{Dynamic Load Balancing:} One area of improvement is the development of dynamic load balancing techniques to distribute work evenly among the available processing units, ensuring that no single unit becomes a bottleneck. This could involve adaptive partitioning of data based on the size and complexity of the tasks, as well as real-time monitoring of the processing load to reassign tasks as needed.
\item \textbf{Advanced Data Partitioning:} Investigating advanced data partitioning strategies that consider the underlying data structures and relationships can lead to more efficient parallel processing. Techniques such as graph partitioning, k-d tree partitioning, or space-filling curves can be explored to optimize the data division process and minimize communication overhead.

\item \textbf{Fault Tolerance and Recovery Mechanisms:} Developing more robust fault tolerance and recovery mechanisms can increase the resilience of the proposed solution in the face of failures, such as node crashes, network issues, or data corruption. Techniques such as checkpointing, replication, and automatic recovery can be explored to ensure the continuous operation of the system.

\item \textbf{The Optimal number of processes:}When using multiprocessing, we observed that the finishing time was significantly reduced. However, we also noticed that the number of scraped products was sometimes smaller than expected compared to a single process. We will further investigated the relationship between the number of processes and the scraping finish time to identify the optimal number of processes
\end{enumerate}

\clearpage
\bibliographystyle{IEEEtran}
\bibliography{mybib}

\begin{thebibliography}{10}
\providecommand{\url}[1]{#1}
\csname url@samestyle\endcsname
\providecommand{\newblock}{\relax}
\providecommand{\bibinfo}[2]{#2}
\providecommand{\BIBentrySTDinterwordspacing}{\spaceskip=0pt\relax}
\providecommand{\BIBentryALTinterwordstretchfactor}{4}
\providecommand{\BIBentryALTinterwordspacing}{\spaceskip=\fontdimen2\font plus
\BIBentryALTinterwordstretchfactor\fontdimen3\font minus \fontdimen4\font\relax}
\providecommand{\BIBforeignlanguage}[2]{{%
\expandafter\ifx\csname l@#1\endcsname\relax
\typeout{** WARNING: IEEEtran.bst: No hyphenation pattern has been}%
\typeout{** loaded for the language `#1'. Using the pattern for}%
\typeout{** the default language instead.}%
\else
\language=\csname l@#1\endcsname
\fi
#2}}
\providecommand{\BIBdecl}{\relax}
\BIBdecl

\bibitem{b1}
\BIBentryALTinterwordspacing
R.~Diouf, E.~N. Sarr, O.~Sall, B.~Birregah, M.~Bousso, and S.~N. Mbaye, ``Web scraping: State-of-the-art and areas of application,'' \emph{2019 IEEE International Conference on Big Data (Big Data)}, pp. 6040--6042, 2019. [Online]. Available: \url{https://api.semanticscholar.org/CorpusID:211297691}
\BIBentrySTDinterwordspacing

\bibitem{b2}
\BIBentryALTinterwordspacing
X.~Deng, P.~Shiralkar, C.~Lockard, B.~Huang, and H.~Sun, ``Dom-lm: Learning generalizable representations for html documents,'' \emph{ArXiv}, vol. abs/2201.10608, 2022. [Online]. Available: \url{https://api.semanticscholar.org/CorpusID:246285527}
\BIBentrySTDinterwordspacing

\bibitem{b3}
\BIBentryALTinterwordspacing
P.~M. Gulhane, A.~Madaan, R.~R. Mehta, J.~Ramamirtham, R.~Rastogi, S.~Satpal, S.~H. Sengamedu, A.~Tengli, and C.~Tiwari, ``Web-scale information extraction with vertex,'' \emph{2011 IEEE 27th International Conference on Data Engineering}, pp. 1209--1220, 2011. [Online]. Available: \url{https://api.semanticscholar.org/CorpusID:13091007}
\BIBentrySTDinterwordspacing

\bibitem{b4}
\BIBentryALTinterwordspacing
B.~Y. Lin, Y.~Sheng, N.~H. Vo, and S.~Tata, ``Freedom: A transferable neural architecture for structured information extraction on web documents,'' \emph{Proceedings of the 26th ACM SIGKDD International Conference on Knowledge Discovery \& Data Mining}, 2020. [Online]. Available: \url{https://api.semanticscholar.org/CorpusID:221191345}
\BIBentrySTDinterwordspacing

\bibitem{b9}
\BIBentryALTinterwordspacing
M.~Dayalan, ``Mapreduce: simplified data processing on large clusters,'' \emph{Commun. ACM}, vol.~51, pp. 107--113, 2008. [Online]. Available: \url{https://api.semanticscholar.org/CorpusID:67055872}
\BIBentrySTDinterwordspacing

\bibitem{b5}
\BIBentryALTinterwordspacing
{Wikipedia contributors}, ``Leafly --- {Wikipedia}{,} the free encyclopedia,'' 2023, [Online; accessed 19-December-2023]. [Online]. Available: \url{https://en.wikipedia.org/w/index.php?title=Leafly&oldid=1157142050}
\BIBentrySTDinterwordspacing

\bibitem{b10}
\BIBentryALTinterwordspacing
M.~A. Zaharia, R.~Xin, P.~Wendell, T.~Das, M.~Armbrust, A.~Dave, X.~Meng, J.~Rosen, S.~Venkataraman, M.~J. Franklin, A.~Ghodsi, J.~E. Gonzalez, S.~Shenker, and I.~Stoica, ``Apache spark: a unified engine for big data processing,'' \emph{Commun. ACM}, vol.~59, pp. 56--65, 2016. [Online]. Available: \url{https://api.semanticscholar.org/CorpusID:251649227}
\BIBentrySTDinterwordspacing

\bibitem{b6}
``Richardson, l. (2007). beautiful soup documentation. april.''

\bibitem{b7}
``Selenium contributors. selenium with python. 2021. https://selenium-python.readthedocs.io/.''

\bibitem{b8}
``Garrett, j. j. (2005). ajax: A new approach to web applications.''

\bibitem{b11}
\BIBentryALTinterwordspacing
M.~Rocklin, ``Dask: Parallel computation with blocked algorithms and task scheduling,'' in \emph{SciPy}, 2015. [Online]. Available: \url{https://api.semanticscholar.org/CorpusID:63554230}
\BIBentrySTDinterwordspacing

\bibitem{b12}
\BIBentryALTinterwordspacing
Z.~Chi, Y.~Li, H.~Sun, Y.~Yao, Z.~Lu, and T.~Zhu, ``B2w2: N-way concurrent communication for iot devices,'' \emph{Proceedings of the 14th ACM Conference on Embedded Network Sensor Systems CD-ROM}, 2016. [Online]. Available: \url{https://api.semanticscholar.org/CorpusID:8208501}
\BIBentrySTDinterwordspacing

\bibitem{b13}
\BIBentryALTinterwordspacing
Z.~Chi, Z.~Huang, Y.~Yao, T.~Xie, H.~Sun, and T.~Zhu, ``Emf: Embedding multiple flows of information in existing traffic for concurrent communication among heterogeneous iot devices,'' \emph{IEEE INFOCOM 2017 - IEEE Conference on Computer Communications}, pp. 1--9, 2017. [Online]. Available: \url{https://api.semanticscholar.org/CorpusID:958378}
\BIBentrySTDinterwordspacing

\bibitem{b14}
\BIBentryALTinterwordspacing
Y.~Li, Z.~Chi, X.~Liu, and T.~Zhu, ``Chiron: Concurrent high throughput communication for iot devices,'' \emph{Proceedings of the 16th Annual International Conference on Mobile Systems, Applications, and Services}, 2018. [Online]. Available: \url{https://api.semanticscholar.org/CorpusID:49668907}
\BIBentrySTDinterwordspacing

\bibitem{80}
W.~Wang, T.~Xie, X.~Liu, and T.~Zhu, ``Ect: Exploiting cross-technology concurrent transmission for reducing packet delivery delay in iot networks,'' in \emph{IEEE INFOCOM 2018-IEEE Conference on Computer Communications}.\hskip 1em plus 0.5em minus 0.4em\relax IEEE, 2018, pp. 369--377.

\bibitem{450}
\BIBentryALTinterwordspacing
Z.~Chi, Y.~Li, Z.~Huang, H.~Sun, and T.~Zhu, ``Simultaneous bi-directional communications and data forwarding using a single zigbee data stream,'' \emph{IEEE INFOCOM 2019 - IEEE Conference on Computer Communications}, pp. 577--585, 2019. [Online]. Available: \url{https://api.semanticscholar.org/CorpusID:86852290}
\BIBentrySTDinterwordspacing

\bibitem{16}
\BIBentryALTinterwordspacing
Z.~Chi, Y.~Yao, T.~Xie, X.~Liu, Z.~Huang, W.~Wang, and T.~Zhu, ``Ear: Exploiting uncontrollable ambient rf signals in heterogeneous networks for gesture recognition,'' \emph{Proceedings of the 16th ACM Conference on Embedded Networked Sensor Systems}, 2018. [Online]. Available: \url{https://api.semanticscholar.org/CorpusID:53092538}
\BIBentrySTDinterwordspacing

\bibitem{44}
\BIBentryALTinterwordspacing
A.~Li, J.~Li, D.~Han, Y.~Zhang, T.~Li, and T.~Zhu, ``Phyauth: Physical-layer message authentication for zigbee networks,'' in \emph{USENIX Security Symposium}, 2023. [Online]. Available: \url{https://api.semanticscholar.org/CorpusID:260340570}
\BIBentrySTDinterwordspacing

\bibitem{22}
\BIBentryALTinterwordspacing
G.~Song and T.~Zhu, ``Ml-based secure low-power communication in adversarial contexts,'' \emph{ArXiv}, vol. abs/2212.13689, 2022. [Online]. Available: \url{https://api.semanticscholar.org/CorpusID:255186145}
\BIBentrySTDinterwordspacing

\bibitem{34}
\BIBentryALTinterwordspacing
W.~Iqbal, W.~Wang, and T.~Zhu, ``Machine learning and artificial intelligence in next-generation wireless network,'' \emph{ArXiv}, vol. abs/2202.01690, 2021. [Online]. Available: \url{https://api.semanticscholar.org/CorpusID:246485682}
\BIBentrySTDinterwordspacing

\bibitem{59}
\BIBentryALTinterwordspacing
J.~Wang, Z.~Huang, W.~Zhang, A.~Patil, K.~Patil, T.~Zhu, E.~J. Shiroma, M.~A. Schepps, and T.~B. Harris, ``2016 ieee international conference on big data (big data) wearable sensor based human posture recognition.'' [Online]. Available: \url{https://api.semanticscholar.org/CorpusID:10607990}
\BIBentrySTDinterwordspacing

\bibitem{58}
\BIBentryALTinterwordspacing
J.~Gao, P.~Yi, Z.~Chi, and T.~Zhu, ``A smart medical system for dynamic closed-loop blood glucose-insulin control,'' \emph{Smart Health}, pp. 18--33, 2017. [Online]. Available: \url{https://api.semanticscholar.org/CorpusID:80388593}
\BIBentrySTDinterwordspacing

\bibitem{52}
\BIBentryALTinterwordspacing
S.~Claros, W.~Wang, and T.~Zhu, ``Investigations of smart health reliability,'' \emph{ArXiv}, vol. abs/2112.15169, 2021. [Online]. Available: \url{https://api.semanticscholar.org/CorpusID:245634369}
\BIBentrySTDinterwordspacing

\bibitem{60}
\BIBentryALTinterwordspacing
S.~Li, P.~Yi, Z.~Huang, T.~Xie, and T.~Zhu, ``Energy scheduling and allocation in electric vehicles energy internet,'' \emph{2016 IEEE Power \& Energy Society Innovative Smart Grid Technologies Conference (ISGT)}, pp. 1--5, 2016. [Online]. Available: \url{https://api.semanticscholar.org/CorpusID:26187112}
\BIBentrySTDinterwordspacing

\bibitem{63}
\BIBentryALTinterwordspacing
Z.~Huang and T.~Zhu, ``eair: an energy efficient air quality management system in residential buildings: poster abstract,'' \emph{Proceedings of the 1st ACM Conference on Embedded Systems for Energy-Efficient Buildings}, 2014. [Online]. Available: \url{https://api.semanticscholar.org/CorpusID:117747}
\BIBentrySTDinterwordspacing

\bibitem{19}
\BIBentryALTinterwordspacing
A.~V. Bhaskar, A.~Baingane, R.~Jahnige, Q.~Zhang, and T.~Zhu, ``A secured protocol for iot networks,'' \emph{ArXiv}, vol. abs/2012.11072, 2020. [Online]. Available: \url{https://api.semanticscholar.org/CorpusID:229339848}
\BIBentrySTDinterwordspacing

\bibitem{68}
\BIBentryALTinterwordspacing
Z.~Zhou, X.~Kuang, L.~Sun, L.~Zhong, and C.~Xu, ``Endogenous security defense against deductive attack: When artificial intelligence meets active defense for online service,'' \emph{IEEE Communications Magazine}, vol.~58, pp. 58--64, 2020. [Online]. Available: \url{https://api.semanticscholar.org/CorpusID:220606123}
\BIBentrySTDinterwordspacing

\bibitem{51}
\BIBentryALTinterwordspacing
Y.~Pan, Q.~Chen, N.~Zhang, Z.~Li, T.~Zhu, and Q.~Han, ``Extending delivery range and decelerating battery aging of logistics uavs using public buses,'' \emph{IEEE Transactions on Mobile Computing}, vol.~22, pp. 5280--5295, 2023. [Online]. Available: \url{https://api.semanticscholar.org/CorpusID:248173721}
\BIBentrySTDinterwordspacing

\bibitem{57}
\BIBentryALTinterwordspacing
Y.~Pan, B.~K. Bhargava, Z.~Ning, N.~Slavov, S.~Li, J.~Liu, S.~Xu, C.~Li, and T.~Zhu, ``Safe and efficient uav navigation near an airport,'' \emph{ICC 2019 - 2019 IEEE International Conference on Communications (ICC)}, pp. 1--6, 2019. [Online]. Available: \url{https://api.semanticscholar.org/CorpusID:198169074}
\BIBentrySTDinterwordspacing

\end{thebibliography}

\end{document}